\documentclass[preprint,12pt]{elsarticle}

\usepackage{amssymb}
\usepackage{amsmath}
\usepackage{amsthm}

\usepackage{tikz}
\usepackage{tikz-3dplot}
\usepackage{pgfplots}
\usepackage{subcaption}
\pgfplotsset{compat=1.6}
\usetikzlibrary{calc,patterns,angles,quotes}
\usepackage{graphicx}
\usepackage{xcolor}
\usepackage{hyperref}
\usepackage{rotating}
\usepackage{makecell} 

\newcommand{\cc}[1]{\overline{#1}}

\theoremstyle{thmstylethree}%
\journal{TODO}

\begin{document}

\begin{frontmatter}



\title{A Generalized Adjoint Method for non-holomorphic cost and constraint functions}


\author[1]{Andrew Zheng} 
\ead{andrewf.zheng@mail.utoronto.ca}
\author[1]{Adam R. Stinchcombe}
\ead{stinch@math.toronto.edu}

\cortext[cor1]{Corresponding author}

\affiliation[1]{organization={Department of Mathematics, University of Toronto},
            addressline={40 St George Street}, 
            city={Toronto},
            postcode={M5S 2E4}, 
            state={Ontario},
            country={Canada}}

\begin{abstract}
The adjoint method is an efficient way to numerically compute gradients in optimization problems with constraints, but is only formulated to differentiable cost and constraint functions on real variables. With the introduction of complex variables, which occur often in many inverse problems in electromagnetism and signal processing problems, both the cost and constraint can become non-holomorphic and hence non-differentiable in the standard definitions. Using the notion of CR-calculus, a generalized adjoint method is introduced that can compute the direction of steepest ascent for the cost function while enforcing the constraint even if both are non-holomorphic. 

\end{abstract}




\end{frontmatter}


\section{Introduction}\label{sec: introduction}

Optimization problems with partial differential equation (PDE) constraints such as inverse problems are complicated to solve numerically. In many gradient based optimization algorithms, common methods to impose the constraints such as Lagrange multipliers and penalty methods are not very effective and can be computationally expensive.
The adjoint method offers an efficient way to compute gradients that enforce the constraint, but at the price that both the cost and constraints are differentiable. 

Suppose we have a cost function $f(x, p)$ we wish to minimize over $p$ with a constraint $g(x, p) = 0$ and we know that for each $p$, there is a unique $x$ such that $g(x, p) = 0$. The adjoint method states that 
\begin{equation}\label{eq: original adjoint method}
    \frac{\partial f}{\partial p} = g_p^* \lambda, \quad g_x^*\lambda = -f_x
\end{equation}
where the subscripts represent partial derivatives with respect to that input, and not a total derivative in case other inputs depend on it.
The benefit of this is that one no longer needs to compute $\frac{\partial x}{\partial p}$ as it is usually very complicated. 
Standard implementations are restricted to functions whose domain and range are both real, where the notation of differentiability is easily defined through vector calculus and limits. However, when dealing with complex variables, our standard definitions no longer work. 

Let $z = x + iy$ and consider a function $f: \mathbb{C} \rightarrow \mathbb{C}$ written as 
\begin{equation}
    f(z) = f(x+iy) = u(x, y) + iv(x, y).
\end{equation}
The derivative of $f$ at a point $z$ can still be defined as 
\begin{equation}
    f'(z) = \lim_{h \rightarrow 0} \frac{f(z+h) - f(z)}{h}.
\end{equation}
However, $h$ can now approach $0$ in multiple ways. 
The standard complex analysis course in mathematics gives the definition that $f$ is holomorphic if it satisfies the Cauchy Riemann equations, which are
\begin{equation}
    \frac{\partial u}{\partial x} = \frac{\partial v}{\partial y}, \quad \frac{\partial v}{\partial x} = -\frac{\partial u}{\partial y}.
\end{equation}
It is well known that $f'(z)$ is well defined and exists if and only if $f$ is holomorphic. 

The problem with this definition is that our objective functions in optimization usually maps to the real line as $\mathbb{C}$ is not well ordered. This means that $v(x, y) = 0$ for all $x, y$. Substituting this into the Cauchy Riemann equations, we get that $f: \mathbb{C} \rightarrow \mathbb{R}$ is holomorphic (and hence differentiable in the standard definition) if and only if it is a constant function. As objective functions are usually non-constant, we need to generalize our gradient based optimization algorithms to non-holomorphic functions. 
A simple example is the function 
\begin{equation}
    f(z) = |z|^2 = z \cc{z}.
\end{equation}
It is not holomorphic, but has an obvious minimum at $z=0$. Furthermore, its graph is smooth and the direction of steepest descent/ascent is intuitively clear. 

Wirtinger calculus was first proposed as an idea on how to define derivatives of non-holomorphic functions, and then formalized into CR-calculus \cite{kreutz2009complex}. However, CR-calculus has yet been fully combined with the adjoint method. 
There have been attempts to allow for complex variables in the adjoint method \cite{plessix2006review, johnson2012notes, virieux2009overview, somersalo1992existence}, but all the constraints under consideration were holomorphic with respect to the complex variable. 

In some considerations of complex variables in the adjoint method, only holomorphic linear constraints are considered \cite{johnson2012notes, plessix2006review}. In this case, the adjoint equations that arises happen to be very similar with the standard adjoint method equations. 
Another common type of optimization problem that is discussed is when $x$ is complex valued, so the cost function is non-holomorphic, but $p$ is real \cite{vogel2002computational, virieux2009overview}. This results in everything being automatically holomorphic in $p$, and the adjoint equations also become very similar to the standard adjoint method equations. 

To our knowledge, there have been little to no discussion of optimization problems where both the cost and constraint functions are non-holomorphic in both $x$ and $p$. 
There are many potential reasons for the lack of research in this area. 
In many PDEs, complex variables are only included as coefficients or boundary conditions, and hence the constraint is holomorphic in these variables. A common example is conductivity in electromagnetics is complex valued, but often treated as real in numerical problems \cite{tyni2024boundary}. 
Complex impedance values appear in geophysics forward problems in a non-linear way \cite{virieux2009overview}, but the model parameters $p$ are assumed to be real, resulting in many functions automatically being holomorphic. 

The second reason is because complex numbers $z = x+iy$ can be viewed as elements $(x, y)\in\mathbb{R}^2$, so it is possible to deal with complex variables by doubling the dimension \cite{johnson2012notes}. This is not recommended due to the following reasons.
Consider a simple linear constraint $g(x, p) = A(p)x -b$. By viewing $\mathbb{C}^n \cong \mathbb{R}^{2n}$, the size of all vectors are doubled, while the size of the matrix is quadrupled. 
Especially for the matrix, the increase in dimension from $n^2$ to $4n^2$ not only increases the memory needed to store a potentially dense matrix, but also can increase the condition number, making solving this linear system very hard.
Furthermore, by storing real and imaginary values in different rows, the matrix has a lot of entries with duplicated values. The matrix 
\begin{equation}
    \begin{bmatrix}
        a + bi & \cdots \\
        \vdots & \ddots
    \end{bmatrix}
\end{equation}
expressed using only real variables would be 
\begin{equation}
    \begin{bmatrix}
    a & -b & \cdots \\
    b & a & \cdots \\
    \vdots & \vdots & \ddots
\end{bmatrix}.
\end{equation}
Sometimes, the spectrum of the original matrix over complex variables has a meaning, but it is lost when representing $\mathbb{C}$ as $\mathbb{R}^2$. 

From a theoretical stand point, though $\mathbb{C}^n$ and $\mathbb{R}^{2n}$ are isomorphic as vector spaces, they are two completely different Hilbert spaces. Manifolds defined on them have very different tangent spaces, which also means their derivatives and gradients cannot be treated as the same.
If the underlying physical phenomenon is defined using the complex field, it would be unnatural to model it using only real variables as the field structures are inherently different.

In this paper, \autoref{sec: foundations} introduces the fundamentals of CR-calculus. \autoref{sec: generalized adjoint method} extends these definitions into the generalized adjoint method that works for non-holomorphic functions using the notation of differential geometry. This was chosen so that the equations can be applied to both the discrete and continuous adjoint method formulations \cite{giannakoglou2008adjoint}. \autoref{sec: lagrangian formulation} gives the same results, but with a Lagrangian formulation of the problem. Finally, we end with some examples and numerical results in \autoref{sec: results}. 

\section{Foundations of complex differentials}\label{sec: foundations}

To really understand how to define the derivatives of non-holomorphic functions, we need to understand the difference between the derivative and the gradient.
In the abstract case, let us consider smooth finite dimensional manifolds $M, N$ and $\mathbb{F}$ a field that is either $\mathbb{C}$ or $\mathbb{R}$. 
For a point $p\in M$, its tangent space $T_pM$ can be viewed as a vector space comprised of tangent vectors of curves passing through the point $p$. 
For any $v\in T_pM$, it can also be viewed as a directional derivative $v: C^{\infty}(M) \rightarrow \mathbb{F}$. 
With the basis derivations $\frac{\partial}{\partial x^k}|_p$, it acts on smooth functions by 
\begin{equation}
    \frac{\partial}{\partial x^k}|_p g = \frac{\partial g}{\partial x^k}(p).
\end{equation}
Let $g: M \rightarrow N$ be a smooth function on $M$ and $p$ a point in $M$. 
The differential of $g$ is a covector field $dg$ defined by 
\begin{equation}
    dg |_p(v) = vg, \quad v\in T_pM
\end{equation}
where $v$ is viewed as a derivation. 

Since $g$ is smooth, it induces a linear map between the tangent spaces $T_p M$ and $T_{g(p)}N$. 
For any curve $\gamma(t) \in M$ with tangent vector $v$ at $p$, the derivative $dg_p$ (or push forward $g_{(*, p)}$) maps $v$ to a tangent vector $dg|_p(v) \in T_{g(p)}N$. 
The function $g$ also induces a linear map between the cotangent spaces $T_{g(p)}^*N$ and $T_p^*M$. 
This map $g^*$ is the pull back that acts on $\alpha\in T_{g(p)}^*N$ via 
\begin{equation}
    g^*(\alpha)(v) = \alpha (g_*v), \quad v\in T_pM.
\end{equation}

To connect this with the chain rule, let $N_2$ be another manifold and let $f: N \rightarrow N_2$ be a smooth function. 
Consider the function $h: M \rightarrow N_2$ be defined by $h = f\circ g$. 
Its differential $dh$ is a pushforward map between $T_p M$ and $T_{f(g(p))}N_2$, but it is also an element in the cotangent space $T_p^*M$.
We have that 
\begin{equation}
    d(f\circ g)_p = (f\circ g)_{(*, p)} = f_{(*, g(p))} \circ g_{(*, p)}.
\end{equation}
Acting on a tangent vector $v\in T_p^N$, we have 
\begin{equation}
    (f\circ g)_{(*, p)}(v) = f_{(*, g(p))}\left(g_{(*, p)}(v)\right) = g^*(f_{(*, g(p))})(v).
\end{equation}
Hence, $d(f\circ g)_p = g^* (df_{g(p)}) = g^*(f_{(*, g(p)})$.

Now that we understand that the differential of a function is a map between tangent spaces, we need to define what the gradient is. 
Given a function $g: M \rightarrow N$, its gradient can be defined as follows. 
Let $\Omega$ be the Riemannian metric of $M$, then the gradient $\nabla g$ is the unique vector field that satisfies 
\begin{equation}
    \langle \nabla g, X\rangle_{\Omega} = Xg
\end{equation}
for every vector field $X$. Equivalently, for each $p\in M$, 
\begin{equation}
    \langle \nabla f|_p, \cdot \rangle_{\Omega} = dg \in T_p^*M. 
\end{equation} 
By Riez representation, the tangent vector $\nabla g(p)$ exists for all $p\in M$ and lies in $T_pM$.

As the definition of the gradient depends on the differential, we now define the differential of a general function $f$ whose domain is $N\subset \mathbb{C}$. Using CR-calculus, we view $z$ and $\cc{z}$ as two independent variables, and define the differential of $f: N \rightarrow \mathbb{F}$ as 
\begin{equation}
    df = \frac{\partial f}{\partial z} dz + \frac{\partial f}{\partial \cc{z}} d\bar{z}. 
\end{equation}
More formally, we define $c = (z, \cc{z})$, and the domain of $f$ is a subspace of $\mathbb{C}^2$ defined by $(N \times\mathbb{C})\cap \{(a, b) \in \mathbb{C}^2 : b = \cc{a}\}$. 
If $g: M \rightarrow N$, the chain rule satisfies 
\begin{equation}
    \frac{\partial (f\circ g)}{\partial z} = \frac{\partial f}{\partial g} \circ \frac{\partial g}{\partial z} + \frac{\partial f}{\partial \cc{g}}\circ \frac{\partial \cc{g}}{\partial z}
\end{equation}
and hence
\begin{equation}
    d(f\circ g) = \left(\frac{\partial f}{\partial g} \circ \frac{\partial g}{\partial z} + \frac{\partial f}{\partial \cc{g}} \circ \frac{\partial \cc{g}}{\partial z}\right)dz + \left(\frac{\partial f}{\partial g} \circ \frac{\partial g}{\partial \cc{z}} + \frac{\partial f}{\partial \cc{g}} \circ \frac{\partial \cc{g}}{\partial \cc{z}}\right)d\bar{z}.
\end{equation}
This is difficult to represent using the notation of pushforwards and pullbacks due to the partial derivatives, so we will use the notation of $\frac{\partial }{\partial z}$ and $\frac{\partial }{\partial \cc{z}}$. 

Returning to the simple case of $f: N \rightarrow \mathbb{F}$, the case of $\mathbb{F} = \mathbb{R}$ gives us that $\frac{\partial f}{\partial \cc{z}} = \overline{\frac{\partial f}{\partial z}}$ and the differential simplifies to 
\begin{equation}
    df = 2 \Re\left( \frac{\partial f}{\partial z} dz \right).
\end{equation}
It acts on a tangent vector $v\in T_pM$ by 
\begin{equation}
    df[v] = 2 \Re \left( \frac{\partial f}{\partial z}[v] \right). 
\end{equation}
Though $df$ acts on tangent vectors, it is not an element of $T_p^*M$. 
Due to taking the real part, it is not complex linear and can only be considered a functional on $T_pM$. 
We can still define the gradient of $f$ at $p$ to be the vector such that 
\begin{equation}
    \Re\left(\langle \nabla f|_p, v\rangle_\Omega \right) = df[v],
\end{equation}
but we can no longer use Riesz representation to give existence and uniqueness. 
There is no general theory that can prove existence and uniqueness for arbitrary Riemannian metrics, and this should not come as a surprise since $f$ is a non-holomorphic function.
In the case where $\langle\cdot, \cdot \rangle_{\Omega}$ is the standard inner product, it is easy to see that 
\begin{equation}\label{eq: gradient with standard inner product}
    \nabla f = 2\cc{\frac{\partial f}{\partial z}}.
\end{equation}

If $\mathbb{F} = \mathbb{C}$ and $f$ is holomorphic, then $\frac{\partial f}{\partial \cc{z}} = 0$ and 
\begin{equation}
    df = \frac{\partial f}{\partial z} dz
\end{equation}
which is what we would expect. 
The gradient in this case is $\nabla f = \cc{\frac{\partial f}{\partial z}}$. 
If $\mathbb{F} = \mathbb{C}$ but $f$ is not holomorphic, we want to define the gradient so that 
\begin{equation}
    \langle \nabla f, v\rangle = df[v].
\end{equation}
This differs from the case when $\mathbb{F} = \mathbb{R}$ as we no longer need the real part since $T_{f(z)}\mathbb{F}$ is a complex Hilbert space. 
For functions defined on $\mathbb{C}^n$, the same ideas work, but the partial derivatives are replaced by Jacobians. 

Consider the simple example of $f(z) = |z|^2$, its differential is 
\begin{equation}
    df[v] = 2 \Re(\cc{z}[v]) = 2\Re (zv)
\end{equation}
and its gradient is easily calculated to be $\nabla f = 2z$ using \autoref{eq: gradient with standard inner product}. 
Comparing this to the holomorphic function $g(z) = z^2$, its differential is $dg = 2z dz$ and its gradient is $2z$. 
Both have a factor of $2$, but for completely different reasons. 
This shows that CR-calculus is not a generalization of complex derivatives for holomorphic functions, but an extension of it to arbitrary functions defined on $\mathbb{C}$. 
Readers are encouraged to read the original notes by Keutz-Delgado \cite{kreutz2009complex} for a more detailed explanation with less differential geometry.


\section{The Generalized Adjoint Method}\label{sec: generalized adjoint method}

Let us consider the cost function $f(x, p)$ with the constraint $g(x, p) = 0$. Since there is a unique $x$ for each $p$ such that $g(x, p) = 0$, we can write $f(x(p), p)$ as $(f\circ F)(p)$ where $F(p) = (x(p), p)$. For notations, we use $\frac{\partial f}{\partial p}$ to denote a partial derivative of the function $f(x(\cdot), \cdot)$ with respect to $p$ and $\partial_1 f$ and $\partial_2 f$ to represent partial derivatives in the first and second arguments of $f(\cdot, \cdot)$. For the partial derivatives of the conjugate variables, we represent them by $\partial_1^c f$ and $\partial_2^c f$. Similarly, we denote $\nabla_k f$ as the gradient of $f$ with respect to the $k^{\text{th}}$ input and $\nabla_k^c f$ as the gradient of $f$ with respect to the conjugate of the $k^{\text{th}}$ input. 

All derivatives are viewed as pushforward maps between tangent spaces of manifolds for the most general formulation of the adjoint method. Hence, all equations can be applied to both continuous or discrete constraints. In the discrete case, partial derivatives are defined using vector calculus. In the continuous case, partial derivatives are defined using Fr\'echet derivatives. 
We now generalize the adjoint method to complex non-holomorphic functions. 

Let $M$ be the manifold in which $p$ belongs to. This can be a complex or a real manifold.
Let $N$ be the manifold in which $x$ belongs to, which is a complex manifold. Both of these manifolds can be of finite dimension (corresponding to the discrete formulation), or infinite dimensions such as function spaces. 
Since $f$ is the objective function, $f$ maps $N\times M$ to $\mathbb{R}$. 
The differential of $f\circ F$ maps tangent vectors in $T_p M$ to $T_{f(x(p), p)} \mathbb{R}$. 
By the chain rule, we have that the differential in $p$ is 
\begin{equation}
    d_p f = \left(\partial_1 f \circ \frac{\partial x}{\partial p} + \partial_1^c f \circ \frac{\partial \bar{x}}{\partial p} + \partial_2 f\right)dp + \left(\partial_1 f \circ \frac{\partial x}{\partial \bar{p}} + \partial_1^c f \circ \frac{\partial \bar{x}}{\partial \bar{p}} + \partial_2^c f\right)d\bar{p}.
\end{equation}
Rearranging the first two terms in the second bracket, we can simplify this to 
\begin{equation}\label{eq: general differential of cost}
     d_p f = 2 \Re \left( \left( \partial_1f \circ \frac{\partial x}{\partial p} + \partial_1^c f \circ \frac{\partial \bar{x}}{\partial p} + \partial_2 f\right)dp\right) 
\end{equation}
The two terms in \autoref{eq: general differential of cost} that are hard to compute are $\frac{\partial x}{\partial p}$ and $\frac{\partial \bar{x}}{\partial p}$. The idea of the adjoint method is to use the constraint to get a formula for these terms that are easier to compute. 

Our constraint $g(x(p), p) = 0$ can also be written as $(g\circ F)(p)=0$. Hence, its differential in $p$ is  
\begin{equation}\label{eq: general differential of constraint}
    d_p g = \left(\partial_1 g \circ \frac{\partial x}{\partial p} + \partial_1^c g \circ \frac{\partial \bar{x}}{\partial p} + \partial_2 g\right)dp + \left(\partial_1 g \circ \frac{\partial x}{\partial \bar{p}} + \partial_1^c g \circ \frac{\partial \bar{x}}{\partial \bar{p}} + \partial_2^c g\right)d\bar{p}.
\end{equation}
As the constraint is forced to be constant, we have that $d_p g = 0$. 
By independence of $dz$ and $d\bar{z}$, this gives us two equations to solve for our two unknowns $\frac{\partial x}{\partial p}$ and $\frac{\partial \bar{x}}{\partial p}$. 
\begin{align}
    0 &= \partial_1 g \circ \frac{\partial x}{\partial p} + \partial_1^c g \circ \frac{\partial \bar{x}}{\partial p} + \partial_2 g, \\
    0 &= \partial_1 g \circ \cc{\frac{\partial \bar{x}}{\partial p}} + \partial_1^c g \circ \cc{\frac{\partial x}{\partial p}} + \partial_2^c g.
\end{align}
Taking the complex conjugate of the second equation and representing the equations using matrices, we have 
\begin{equation}
    -\begin{bmatrix}
        \partial_2 g\\[0.2em]
        \cc{\partial_2^c g}
    \end{bmatrix} = \begin{bmatrix}
        \partial_1g & \partial_1^c g \\[0.2em]
        \cc{\partial_1^c g} & \cc{\partial_1 g} 
    \end{bmatrix} \begin{bmatrix}
        \frac{\partial x}{\partial p}\\[0.2em]
        \frac{\partial \bar{x}}{\partial p}
    \end{bmatrix}
\end{equation}
Hence, 
\begin{equation}\label{eq: generalized adjoint equation}
    \begin{bmatrix}
        \frac{\partial x}{\partial p}\\[0.2em]
        \frac{\partial \bar{x}}{\partial p}
    \end{bmatrix} = \begin{bmatrix}
        \partial_1g & \partial_1^c g \\[0.2em]
        \cc{\partial_1^c g} & \cc{\partial_1 g} 
    \end{bmatrix}^{-1} \left(-\begin{bmatrix}
        \partial_2 g\\[0.2em]
        \cc{\partial_2^c g}
    \end{bmatrix}\right).
\end{equation}
The solution can then be substituted into \autoref{eq: general differential of cost}. 
Just like the original adjoint method, we need to assume that the matrix in \autoref{eq: generalized adjoint equation} is invertible. 
Let $p_0\in M$ and $v\in T_{p_0}M$.
By $\mathbb{R}$-linearity, we can simplify \autoref{eq: general differential of cost} applied to $v$ to 
\begin{equation}
    d_pf[v] = 2\Re\left(\partial_1f \circ \frac{\partial x}{\partial p}[v]\right) + 2\Re\left(\partial_1^c f \circ \frac{\partial \bar{x}}{\partial p}[v]\right) + 2\Re \left(\partial_2 f[v]\right).
\end{equation}
Assuming that we have already computed $\nabla_1 f, \nabla_1^c f$ and $\nabla_2 f$, we get 
\begin{equation}
\begin{split}
    d_pf[v] &= \Re\left(\left\langle \nabla_1 f, \frac{\partial x}{\partial p}[v]\right\rangle\right) +\Re\left(\left\langle \frac{1}{2}\nabla_1^c f, \frac{\partial \bar{x}}{\partial p}[v]\right\rangle\right) + \Re\left(\left\langle\nabla_2 f, v\right\rangle\right)\\
    &= \Re\left(\left\langle \left(\frac{\partial x}{\partial p}\right)^* \nabla_1 f + \left(\frac{\partial \bar{x}}{\partial p}\right)^* \nabla_1^c f + \nabla_2 f, v \right\rangle\right).
\end{split}
\end{equation}
Since $f$ maps to $\mathbb{R}$ and the Riemannian metric is just the standard inner product on $\mathbb{C}$, \autoref{eq: gradient with standard inner product} gives that that the gradient corresponding to $\partial_1^c f$ is the complex conjugate of the gradient corresponding to $\partial_1 f$.
Thus, we can define the gradient of $f$ with respect to $p$ to be 
\begin{equation}\label{eq: gradient for most general case}
    \nabla_p f = \left(\frac{\partial x}{\partial p}\right)^* \nabla_1 f + \left(\frac{\partial \bar{x}}{\partial p}\right)^* \cc{\nabla_1 f} + \nabla_2 f.
\end{equation}
If we fix $\|v\| = 1$, it can easily be seen that the maximum change in $f$ is obtained when we choose $v$ to be the normalized gradient.
Given a point $p_0 \in M$, it is also important to note that $\nabla_p f|_{p_0}$ must be an element in $T_{p_0}M$. If $T_{p_0}M$ is a real Hilbert space, then $\nabla_p f|_{p_0}$ must also be real. This can be enforced by using the fact that $v\in T_{p_0}M$ is also real, so 
\begin{align*}
    d_p f[v] &=  \Re \left(\left\langle \left(\frac{\partial x}{\partial p}\right)^* \nabla_1 f + \left(\frac{\partial \bar{x}}{\partial p}\right)^* \nabla_1^c f + \nabla_2 f, v \right\rangle\right)\\
    &=  \left\langle \Re\left(\left(\frac{\partial x}{\partial p}\right)^* \nabla_1 f + \left(\frac{\partial \bar{x}}{\partial p}\right)^* \nabla_1^c f + \nabla_2 f\right), v \right\rangle. 
\end{align*}
Hence, in the case where $T_{p_0}M$ is a real Hilbert space, 
\begin{equation}
    \nabla_p f = \Re\left(\left(\frac{\partial x}{\partial p}\right)^* \nabla_1 f + \left(\frac{\partial \bar{x}}{\partial p}\right)^* \cc{\nabla_1 f} + \nabla_2 f\right).
\end{equation}
An example that elucidates this will be presented in \autoref{sec: linear constraint example}.

\subsection{Holomorphic in $x$ but not $p$}
Consider the case where the constraint $g(x, p)$ is holomorphic in $x$ but not in $p$, and $x(p)$ is also not holomorphic in $p$. Examples are when $p$ is complex valued, and $g(x,p)$ is a linear system in $x$ but not in $p$. 
We have that \autoref{eq: general differential of constraint} simplifies to 
\begin{equation}
    d_p g = \left(\partial_1 g \circ \frac{\partial x}{\partial p} + \partial_2 g\right)dp + \left(\partial_1 g \circ \frac{\partial x}{\partial \bar{p}} + \partial_2^c g\right)d\bar{p}.
\end{equation}
Setting the differential to zero, we get linear system for $\frac{\partial x}{\partial p}$ and $\frac{\partial x}{\partial \bar{p}}$, but it is much simpler than \autoref{eq: generalized adjoint equation}. 
Using matrices, we have 
\begin{equation}\label{eq: adjoint for holomorphic in x but not p}
    \begin{bmatrix}
        \frac{\partial x}{\partial p}\\[0.2em]
        \frac{\partial \bar{x}}{\partial p}
    \end{bmatrix} = \begin{bmatrix}
        \partial_1g & 0 \\
        0 & \cc{\partial_1 g} 
    \end{bmatrix}^{-1} \left(-\begin{bmatrix}
        \partial_2 g\\[0.2em]
        \cc{\partial_2^c g}
    \end{bmatrix}\right).
\end{equation}
The final gradient still uses \autoref{eq: gradient for most general case}.

\subsection{When everything but the cost function is holomorphic}
If we instead had a holomorphic constraint, meaning that $g$ is holomorphic in $x$ and $p$, and we also assume that $x$ is holomorphic in $p$, then \autoref{eq: generalized adjoint equation} simplifies to 
\begin{equation}
    \begin{bmatrix}
        \frac{\partial x}{\partial p}\\
        0
    \end{bmatrix} = \begin{bmatrix}
        \partial_1g & 0 \\
        0 & \cc{\partial_1 g} 
    \end{bmatrix}^{-1} \left(-\begin{bmatrix}
        \partial_2 g\\
       0
    \end{bmatrix}\right).
\end{equation}
Basically, the second equation becomes trivially true and we are left with 
\begin{equation}\label{eq: adjoint for everything is holomorphic}
    \frac{\partial x}{\partial p} = (\partial_1 g)^{-1}\circ (-\partial_2 g).
\end{equation}
The differential of $f$ also simplifies into 
\begin{equation}\label{eq: d_p f with holomorphic constraint}
    d_p f = 2 \Re \left( \left( \partial_1f \circ \frac{\partial x}{\partial p} + \partial_2 f\right)dp\right)
\end{equation}
Combining \autoref{eq: adjoint for everything is holomorphic} and \autoref{eq: d_p f with holomorphic constraint}, we get 
\begin{equation}\label{eq: d_p f simplified version with holomorphic constraint}
    d_p f = 2 \Re\left(\left(\partial_1 f \circ (\partial_1 g)^{-1}\circ (-\partial_2 g) + \partial_2 f\right)dp\right).
\end{equation}
Given a tangent vector $v \in T_{p_0}M$ for a point $p_0\in M$, we have 
\begin{equation}
    \begin{split}
        d_pf[v] &= 2 \Re\left(\left(\partial_1 f \circ (\partial_1 g)^{-1}\circ (-\partial_2 g) + \partial_2 f\right)[v]\right)\\
        &= 2\Re\left(\partial_1 f \left[(\partial_1 g)^{-1}\circ (-\partial_2 g) [v]\right] \right) + 2 \Re (\partial_2 f [v]).
    \end{split}
\end{equation}
This result is very similar to the standard adjoint method. 

\subsection{Holomorphic linear constraints}\label{sec: linear constraint example}
Many PDE constraints can be discretized into linear constraints, so many inverse problems have constraints of the form $g(x, p) = A(p)x - b$. In this case, $g$ is clearly holomorphic in $x$, but we now also assume that $x$ is holomorphic in $p$. Common examples are if $p$ is real, so everything is automatically holomorphic in $p$, or if $g(x, p)$ is also linear in $p$. This results in
\begin{equation}
    \partial_1 g = A(p), \quad \partial_2 g = \frac{\partial A(p)x}{\partial p}. 
\end{equation}
Using \autoref{eq: d_p f simplified version with holomorphic constraint}, 
\begin{equation}
    \frac{\partial }{\partial p} f [v] = 2\Re\left(\partial_1 f \left[\left(A(p)\right)^{-1}\circ (-\frac{\partial A(p)x}{\partial p}) [v]\right] \right) + 2 \Re (\partial_2 f v)
\end{equation}
We can replace $\partial_1 f$ with $\nabla_1 f$ and the inner product to get 
\begin{equation}\label{eq: differential of f with linear constraint}
    \begin{split}
        d_p f [v] &= \Re\left(\left\langle \nabla_1 f, \left(A(p)\right)^{-1}\circ (-\frac{\partial A(p)x}{\partial p}) [v]\right\rangle \right) + 2 \Re (\partial_2 f v)\\
        &= \Re\left(\left\langle (A(p)^*)^{-1}\nabla_1 f, -\frac{\partial A(p)x}{\partial p}[v]\right\rangle \right) + \Re (\langle \nabla_2 f, v\rangle)\\ 
        &= \Re\left(\left\langle \left(-\frac{\partial A(p)x}{\partial p}\right)^* (A(p)^*)^{-1}\nabla_1 f  + \nabla_2 f, v\right\rangle \right)
    \end{split}
\end{equation}
If $N$ is a real Hilbert space, then the gradient of $f$ with respect to $p$ is 
\begin{equation}
    \nabla_p f = \Re \left(\left(-\frac{\partial A(p)x}{\partial p}\right)^* (A(p)^*)^{-1}\nabla_1 f  + \nabla_2 f\right).
\end{equation}
If $N$ is a complex Hilbert space, then 
\begin{equation}
    \nabla_p f = \left(\left(-\frac{\partial A(p)x}{\partial p}\right)^* (A(p)^*)^{-1}\nabla_1 f  + \nabla_2 f\right).
\end{equation}
In both cases, the adjoint problem is $(A(p)^*)^{-1}\nabla_1 f$ and we need to find a $\lambda$ satisfying 
\begin{equation}
    A(p)^* \lambda = \nabla_1 f.
\end{equation}


\section{A Lagrangian Formulation}\label{sec: lagrangian formulation}

It is also possible to view the adjoint method from a Lagrangian formulation. 
Let $M, N, Q$ be manifolds over $\mathbb{R}$ or $\mathbb{C}$ with $x \in N, p\in M$ and $g(x, p) \in Q$.
Given our functions $f: N\times M \rightarrow \mathbb{R}$ and $g: N\times M \rightarrow Q$, we can define the Lagrangian 
\begin{equation}
    L(x, p, \lambda_1, \lambda_2) = f(x, p) - \lambda_1[g(x, p)] - \lambda_2[\cc{g(x, p)}]
\end{equation}
where $\lambda_1, \lambda_2 \in T^*Q$. It is common to associate the $\lambda_k$ terms as Lagrange multipliers, but technically they are not \cite{givoli2021tutorial}. This is because later on, we do not solve for them along with $x$ and $p$, but we introduce an artificial condition to simplify certain terms in the differentials. This is a Lagrangian formulation, but we are not using Lagrange multipliers.

Since $\lambda_1, \lambda_2$ are elements in the dual space, Riesz Representation theorem states that there exists a $\Tilde{\lambda}_1$ and $\Tilde{\lambda}_2$ such that 
\begin{equation}
    \lambda_1[v] = \langle v, \Tilde{\lambda}_1\rangle_Q, \quad \lambda_2[v] = \langle v, \Tilde{\lambda}_2\rangle_Q
\end{equation}
where $\langle \cdot, \cdot \rangle_Q$ is the inner product defined in $Q$. 

We need to introduce two $\lambda_k$ variables for this single constraint because the Lagrangian needs to map to $\mathbb{R}$. If we only had the single term of $\Re(\lambda[g(x, p)])$, we cannot enforce that the imaginary part of $g(x, p)$ also equals to zero. 
Another way to view this is through the definition of the dual space. Because we are viewing $Q$ as $Q\times \mathbb{C}^{\text{dim}(Q)}\cap \{(a, b): a = \cc{b}\}$, $\lambda$ being in the dual space implies that $\lambda$ lives in a $\text{dim}(Q)$ subspace of a $\text{dim}(2Q)$ space and it can be represented by $\lambda = (\lambda_1,\lambda_2)$ with $\lambda_2 = \cc{\lambda_1}$. Here the dimension of $Q$ can be infinite, as $Q$ can represent some function space. 

Let $h$ be a direction in which we change $\lambda_1$, then 
\begin{equation}
    L(x, p, \lambda_1 + \epsilon h, \lambda_2) - L(x, p, \lambda_1, \lambda_2) = 
    - \langle g(x, p), \epsilon h\rangle_Q.
\end{equation}
Using the definition of Fr\'echet derivatives, the differential in $\lambda_1$ is
\begin{equation}
    d_{\lambda_1} L = -\langle g(x, p), 1 \rangle_Q d\lambda_1 
\end{equation}
For this to be zero, we must have that $g(x, p) = 0$, so equilibrium points of the Lagrangian enforce the constraint. The differential in $\lambda_2$ enforces that $\cc{g(x, p)} = \cc{0}$, so it is the same constraint if we take the complex conjugate. 

As all following inner products are in $M$, we drop the subscript $M$ for ease of notation.
For $v \in T_p M$, we have that 
\begin{equation}
    \begin{split}
        \frac{\partial }{\partial p} \lambda_1 [g(x, p)] [v]&= \left\langle \Tilde{\lambda}_1, \frac{\partial g(x, p)}{\partial p}[v]\right\rangle\\
        &= \left\langle \left(\frac{\partial g(x, p)}{\partial p}\right)^* \Tilde{\lambda}_1, v\right\rangle\\
        &= \left\langle \left(\partial_1 g \circ \frac{\partial x}{\partial p} + \partial_2 g\right)^* \Tilde{\lambda}_1, v\right\rangle.
    \end{split}
\end{equation}
Hence, the differential in $p$ can be written as 
\begin{equation}
\begin{split}
    d_p L &= \bigg(\partial_1 f \frac{\partial x}{\partial p} + \partial_1^c f \circ \frac{\partial \bar{x}}{\partial p} + \partial_2 f \\
    &- \left\langle \left(\partial_1 g \circ \frac{\partial x}{\partial p} + \partial_1^c g \circ \frac{\partial \bar{x}}{\partial p} + \partial_2 g\right)^* \Tilde{\lambda}_1, \cdot\right\rangle\\
    &- \left\langle \left(\partial_1 \bar{g} \circ \frac{\partial x}{\partial p} + \partial_1^c \bar{g} \circ \frac{\partial \bar{x}}{\partial p} + \partial_2 \bar{g}\right)^* \Tilde{\lambda}_2, \cdot\right\rangle \bigg) dp \\
    &+ \bigg(\partial_1 f \frac{\partial x}{\partial \bar{p}} + \partial_1^c f \circ \frac{\partial \bar{x}}{\partial \bar{p}} + \partial_2^c f \\
    &- \left\langle \left(\partial_1 g \circ \frac{\partial x}{\partial \bar{p}} + \partial_1^c g \circ \frac{\partial \bar{x}}{\partial \bar{p}} + \partial_2^c g\right)^* \Tilde{\lambda}_1, \cdot\right\rangle\\
    &- \left\langle \left(\partial_1 \bar{g} \circ \frac{\partial x}{\partial \bar{p}} + \partial_1^c \bar{g} \circ \frac{\partial \bar{x}}{\partial \bar{p}} + \partial_2^c \bar{g}\right)^* \Tilde{\lambda}_2, \cdot\right\rangle \bigg) d\bar{p}. 
\end{split}
\end{equation}
The two terms are not necessarily complex conjugates of each other due to the $\Tilde{\lambda}_k$ terms, so we cannot combine them. 
Assuming we have already calculated $\nabla_1 f, \nabla_1^c f$ and $\nabla_2 f$, we simplify this to 
\begin{equation}
    \begin{split}
        d_p L &= \bigg\langle 
        \left(\frac{\partial x}{\partial p}\right)^*\circ 
        \left(\nabla_1 f - \left(\partial_1 g\right)^*\Tilde{\lambda}_1 - \left(\cc{\partial_1^c g}\right)^*\Tilde{\lambda}_2\right) \\
        &+ \left(\frac{\partial \bar{x}}{\partial p}\right)^*\circ 
        \left(\nabla_1^c f - \left(\partial_1^c g\right)^*\Tilde{\lambda}_1 - \left(\cc{\partial_1 g}\right)^*\Tilde{\lambda}_2\right)\\
        &- \left(\partial_2 g\right)^*\Tilde{\lambda}_1 - \left(\cc{\partial_2^c g}\right)^*\Tilde{\lambda}_2+ \nabla_2 f, \cdot \bigg\rangle dp \\
       &+ \bigg\langle 
        \left(\frac{\partial x}{\partial \bar{p}}\right)^*\circ 
        \left(\nabla_1 f - \left(\partial_1 g\right)^*\Tilde{\lambda}_1 - \left(\cc{\partial_1^c g}\right)^*\Tilde{\lambda}_2\right) \\
        &+ \left(\frac{\partial \bar{x}}{\partial \bar{p}}\right)^*\circ 
        \left(\nabla_1^c f - \left(\partial_1^c g\right)^*\Tilde{\lambda}_1 - \left(\cc{\partial_1 g}\right)^*\Tilde{\lambda}_2\right)\\ 
        &- \left(\partial_2^c g\right)^*\Tilde{\lambda}_1 - \left(\cc{\partial_2 g}\right)^*\Tilde{\lambda}_2 + \nabla_2^c f, \cdot \bigg\rangle d\bar{p}.
    \end{split}
\end{equation}
Here, we also used the fact that $\partial_1 \bar{g} = \cc{\partial_1^c g}$ and $\partial_1^c \bar{g} = \cc{\partial_1 g}$. 

The idea in the Lagrange formulation is to make $d_p L$ not depend on the complicated partials of $x$ and $\bar{x}$ with respect to $p$ and $\bar{p}$. Hence we impose the following equations
\begin{align}
    0 &= \nabla_1 f - \left(\partial_1 g\right)^* \Tilde{\lambda}_1 - \left(\cc{\partial_1^c g}\right)^*\Tilde{\lambda}_2\\
    0 &= \nabla_1^c f - \left(\partial_1^c g\right)^* \Tilde{\lambda}_1 - \left(\cc{\partial_1 g}\right)^* \Tilde{\lambda}_2.
\end{align}
In matrix form, this can be written as 
\begin{equation}
    0 = \begin{bmatrix}
        \left(\partial_1 g\right)^* & \left(\cc{\partial_1^c g}\right)^* \\[0.2em]
        \left(\partial_1^c g\right)^* & \left(\cc{\partial_1 g}\right)^*
    \end{bmatrix} \begin{bmatrix}
        \Tilde{\lambda}_1 \\
        \Tilde{\lambda}_2
    \end{bmatrix} = \begin{bmatrix}
        \nabla_1 f\\
        \nabla_1^c f
    \end{bmatrix}.
\end{equation}
Setting $\Tilde{\lambda}_1$ and $\Tilde{\lambda}_2$ to the solution of this system of equations, we get 
\begin{equation}
\begin{split}
    d_p L &= \Re\left(\left\langle \left(-\begin{bmatrix}
        \partial_2 g\\[0.2em]
        \cc{\partial_2^c g}
    \end{bmatrix}\right)^*\begin{bmatrix}
        \left(\partial_1 g\right)^* & \left(\cc{\partial_1^c g}\right)^* \\[0.2em]
        \left(\partial_1^c g\right)^* & \left(\cc{\partial_1 g}\right)^*
    \end{bmatrix}^{-1}\begin{bmatrix}
        \nabla_1 f\\
        \nabla_1^c f
    \end{bmatrix}+ \nabla_2 f , \cdot \right\rangle dp\right).
\end{split}
\end{equation}
Since $\Tilde{\lambda}_1 = \cc{\Tilde{\lambda}_2}$, the two terms are complex conjugates of each other. If we combine \autoref{eq: generalized adjoint equation} and \autoref{eq: gradient for most general case}, we get 
\begin{equation}
    \nabla_p f = \left(\left(-\begin{bmatrix}
        \partial_2 g\\[0.2em]
        \cc{\partial_2^c g}
    \end{bmatrix}\right)^* \left( \begin{bmatrix}
        \partial_1g & \partial_1^c g \\[0.2em]
        \cc{\partial_1^c g} & \cc{\partial_1 g} 
    \end{bmatrix}^{-1}\right)^* \begin{bmatrix}
        \nabla_1 f\\
        \nabla_1^c f
    \end{bmatrix} + \nabla_2 f\right).
\end{equation}
This agrees with our Lagrange formulation of the generalized adjoint method.
The Lagrange formulation for the case of holomorphic constraints is similar, but also easier as we no longer need to include a $\lambda_2 = \cc{\lambda_1}$ term.


\section{Examples}\label{sec: results}

We now present some numerical examples of computing the gradient with the generalized adjoint method.

\subsection{Holomorphic Example}
Let us consider the constrained minimization problem 
\begin{equation}
    \min f(z, p) = \min |z|^2
\end{equation}
with the constraint 
\begin{equation}\label{eq: linear holomorphic constraint}
    g(z, p) = \begin{bmatrix}
        1 - p_2^2 & 5p_1^2 - 2p_2^2 & 4(p_2 - p_1)\\
        0 & 1-0.1p_1^2 & -50p_2^2\\
        0.1p_1p_2 & p_2^2 + p_1^2 & 1 - 0.75(p_1 + p_2)
    \end{bmatrix}z - \begin{bmatrix}
        0 \\
        0.5\\
        0.5 - 0.5i
    \end{bmatrix} = 0.
\end{equation}
Fixing our domain of $p$ to be $[-0.5, 0.5]\times [-0.5, 0.5]\subset \mathbb{R}^2$, so that the matrix is invertible, we have that $z(p)\in \mathbb{C}^3$ due to the complex vector in $g$. Letting the matrix be denoted as $A(p)$, we have that 
\begin{equation}
    \frac{\partial A(p)z}{\partial p} = \begin{bmatrix}
        10p_1 z_2 - 4z_3 & 2p_2z_1 - 4p_2z_2 + 4z_3\\
        -0.2p_1z_2 & -100p_2z_3\\
        0.1p_2z_1 + 2p_1z_2 - 0.75z_3 & 0.1p_1z_1 + 2p_2z_2 - 0.75z_3
    \end{bmatrix}.
\end{equation}
We also have that $\frac{\partial f}{\partial z} = \Bar{z}$, and it is easy to see that $\nabla_z f = 2z$ satisfies our condition for the gradient. 
Using \autoref{eq: differential of f with linear constraint}, we can compute $d_p f$ at different values of $p$. Since $p$ is on a real manifold, the gradient is computed at the point $p$ by first solving $z = A(p)^{-1}b$, and then using 
\begin{equation*}
    \nabla_p f = \Re \left(\left(-\frac{\partial A(p)z}{\partial p}\right)^* (A(p)^*)^{-1}2z\right).
\end{equation*} 

\begin{figure}
    \centering
    \includegraphics[width=0.75\linewidth]{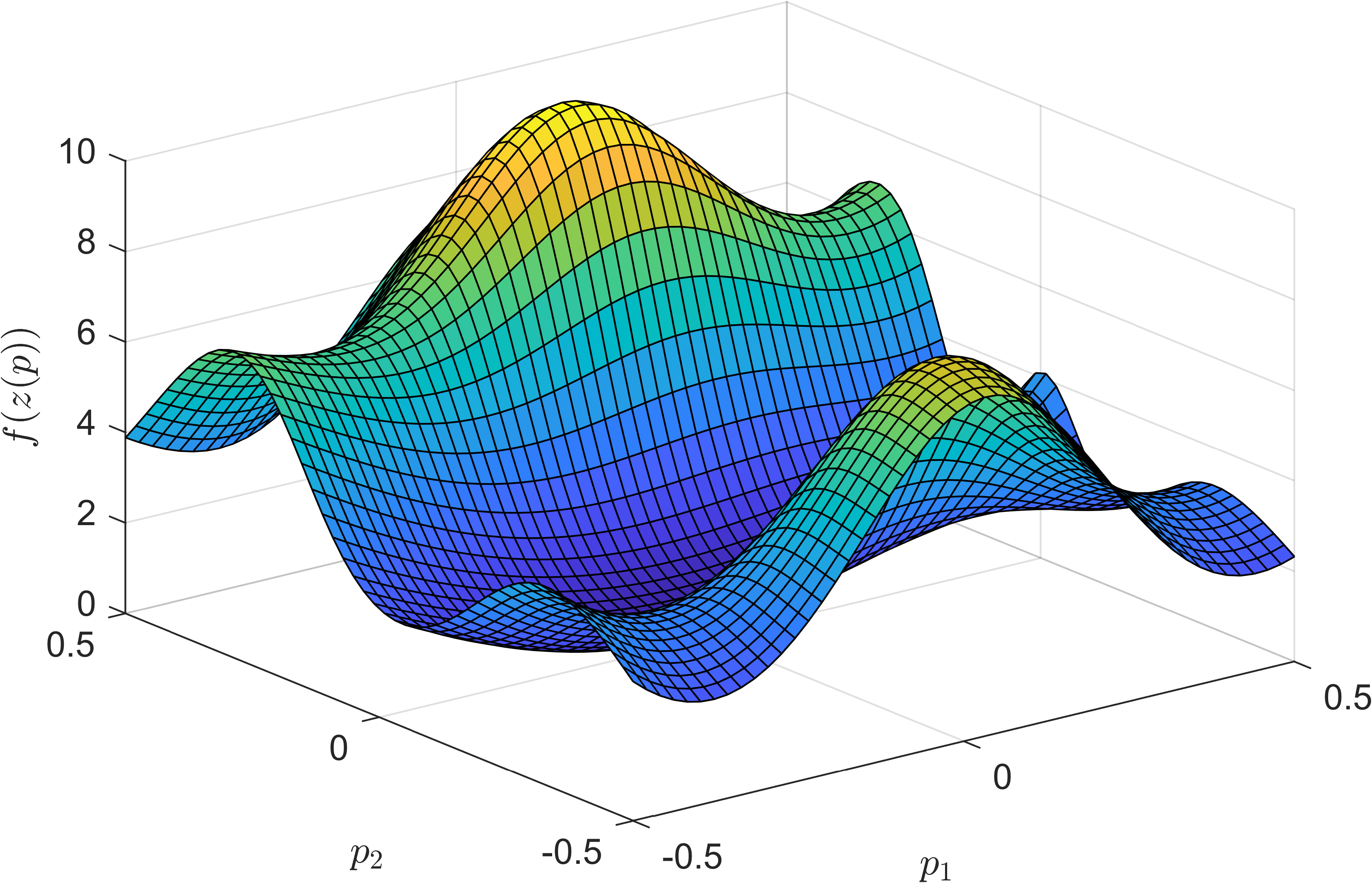}
    \caption{A surface plot of the cost function $f(z(p)) = |z|^2$ with constraint \autoref{eq: linear holomorphic constraint}.}
    \label{fig: surf plot}
\end{figure}

\begin{sidewaysfigure}
    \centering
    \includegraphics[width=0.8\textheight]{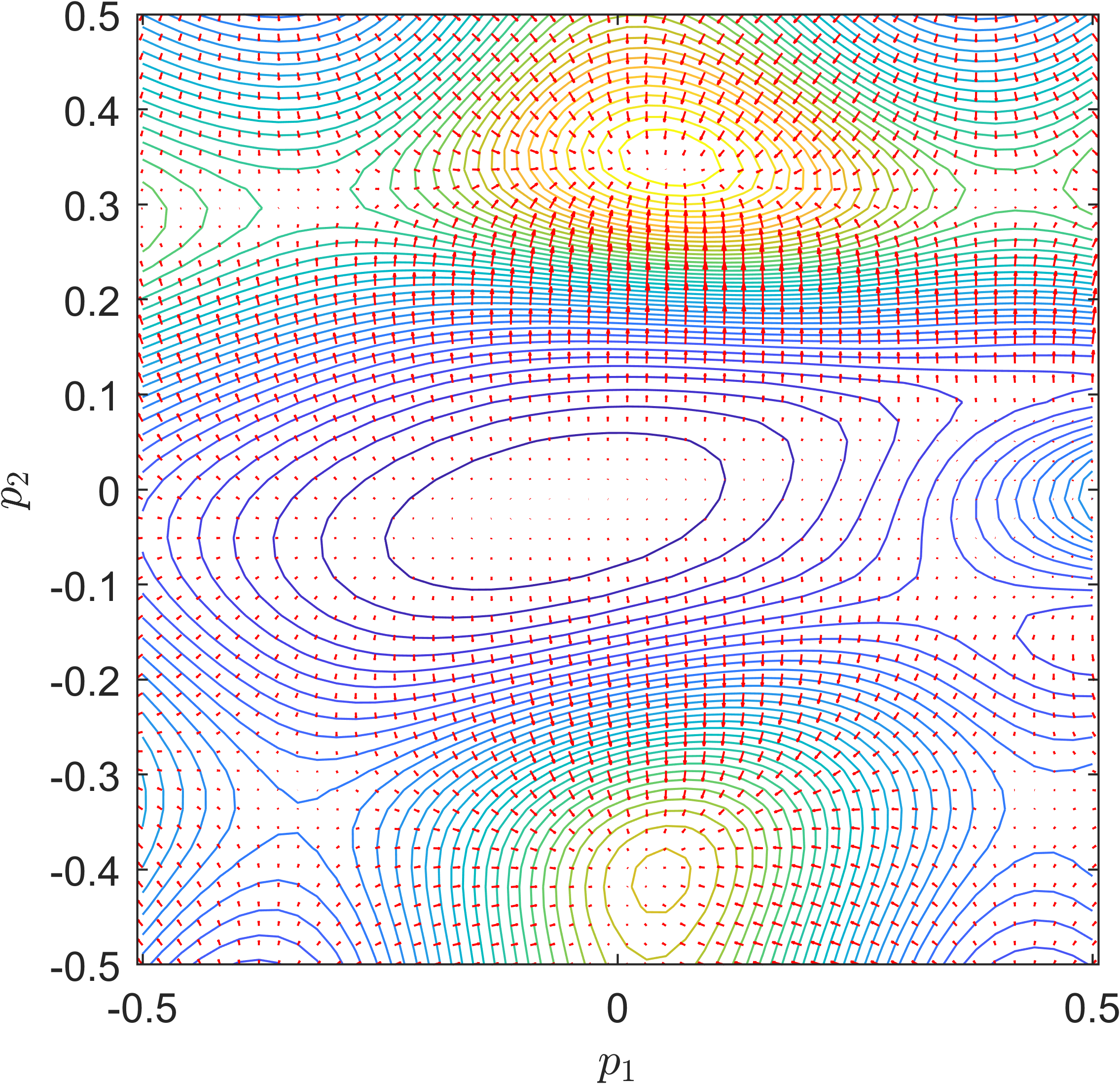}
    \caption{A contour plot of the cost function $f(z(p)) = |z|^2$ with constraint \autoref{eq: linear holomorphic constraint}. The gradients are orthogonal to the contour lines, and points in the direction of steepest ascent. The magnitude of the gradients are proportional to each other, but scaled down so the vectors do not overlap.}
    \label{fig: contour plot}
\end{sidewaysfigure}

As shown in \autoref{fig: contour plot}, the computed gradients using the adjoint method are orthogonal to the contour lines, and we can easily see the location of saddles and local minimum/maximums. Comparing to \autoref{fig: surf plot}, the location vectors are pointing in the correct directions. 
As a function of $p$, $f(z(p))$ is a map from $\mathbb{R}^2 \rightarrow \mathbb{R}$, so we can also compute its gradient using symbolic differentiation. By explicitly solving for $z(p)$ and substituting it into $|z|^2$, we compute a single gradient without any chain rule or adjoint method, and its results match our adjoint method results with machine precision ($\approx 10^{-16}$). We note that for actual complex non-holomorphic functions, some symbolic differentiation packages are wrong. An example is with Matlab's symbolic differentiation of $|z|^2$. It gives $z + \Bar{z}$, which is now a real vector and clearly not the direction of steepest ascent. 

\subsection{Non-holomorphic Example}

Now lets consider a very non-holomorphic example. We define our constraint to be very similar as before, but now we let $p$ be a complex number and define
\begin{equation}
    A(p) = \begin{bmatrix}
        1 - \left(\frac{p - \Bar{p}}{2i}\right)^2 & 5\left(\frac{p + \Bar{p}}{2}\right)^2 - 2\left(\frac{p - \Bar{p}}{2i}\right)^2 & 4 \left(\frac{p - \Bar{p}}{2i} - \frac{p + \Bar{p}}{2}\right)\\
        0 & 1 - 0.1\left(\frac{p + \Bar{p}}{2}\right)^2 & -50 \left(\frac{p - \Bar{p}}{2i}\right)^2\\
        0.1\left(\frac{p - \Bar{p}}{2i} \frac{p + \Bar{p}}{2}\right) & p\Bar{p} & 1 - 0.75\left(\frac{p - \Bar{p}}{2i} + \frac{p + \Bar{p}}{2}\right)
    \end{bmatrix}, \quad b = \begin{bmatrix}
        0\\
        0.5\\
        0.5
    \end{bmatrix}.
\end{equation}
Let $[v_1, v_2, v_3]^T = A(p)^{-1}b \in \mathbb{R}^3$, then we define the cost function $f$ to be 
\begin{equation}
    f(v_1, v_2, v_3) = \|v_1 + v_2i\|^2 + v_3^2.
\end{equation}
If we let $z\in \mathbb{C}^2$, the constraint can be written as 
\begin{equation}\label{eq: non holomorphic constraint}
    g(z, p) = A(p) \begin{bmatrix}
        \frac{z_1 + \Bar{z}_1}{2}\\
        \frac{z_1 - \Bar{z}_1}{2i}\\
        z_2
    \end{bmatrix} - b = 0
\end{equation}
and the cost function is just $f(z) = \|z_1\|^2 + \|z_2\|^2$.
We now compute all required partial derivatives. 
\begin{equation}
    \partial_1 g =  \begin{bmatrix}
        \frac{1}{2}\left(1 - \left(\frac{p - \Bar{p}}{2i}\right)^2\right) + \frac{1}{2i}\left( 5\left(\frac{p + \Bar{p}}{2}\right)^2 - 2\left(\frac{p - \Bar{p}}{2i}\right)^2\right) & 4 \left(\frac{p - \Bar{p}}{2i} - \frac{p + \Bar{p}}{2}\right)\\
        \frac{1}{2i}\left( 1 - 0.1\left(\frac{p + \Bar{p}}{2}\right)^2 \right) & -50 \left(\frac{p - \Bar{p}}{2i}\right)^2\\
        \frac{1}{2} \left(0.1\left(\frac{p - \Bar{p}}{2i} \frac{p + \Bar{p}}{2}\right)\right) + \frac{1}{2i}\left( p\Bar{p}\right) & 1 - 0.75\left(\frac{p - \Bar{p}}{2i} + \frac{p + \Bar{p}}{2}\right)        
    \end{bmatrix}.
\end{equation}
\begin{equation}
    \partial_1^c g = \begin{bmatrix}
        \frac{1}{2}\left(1 - \left(\frac{p - \Bar{p}}{2i}\right)^2\right) - \frac{1}{2i}\left( 5\left(\frac{p + \Bar{p}}{2}\right)^2 - 2\left(\frac{p - \Bar{p}}{2i}\right)^2\right) & 0 \\
        -\frac{1}{2i}\left( 1 - 0.1\left(\frac{p + \Bar{p}}{2}\right)^2 \right) & 0 \\
        \frac{1}{2} \left(0.1\left(\frac{p - \Bar{p}}{2i} \frac{p + \Bar{p}}{2}\right)\right) - \frac{1}{2i}\left( p\Bar{p}\right) & 0        
    \end{bmatrix}.
\end{equation}
Note that $\partial_1 g \ne \overline{\partial_1^c g}$ as $g$ does not map to real values. Though $A(p)^{-1}b$ gives real values, the general function $g(z, p)$ does not map to real values. 

\begin{equation}
    \partial_2 g = \begin{bmatrix}
        -\frac{2}{2i}\frac{p - \Bar{p}}{2i}\frac{z_1 + \Bar{z}_1}{2} + \left(\frac{10}{2}\frac{p+\Bar{p}}{2} - \frac{4}{2i}\frac{p-\Bar{p}}{2i}\right)\frac{z_1 - \Bar{z}_1}{2i} + 4\left(\frac{1}{2i} - \frac{1}{2}\right)z_2\\
        -\frac{0.2}{2}\frac{p+\Bar{p}}{2}\frac{z_1 - \Bar{z}_1}{2i} - \frac{100}{2i}\frac{p-\Bar{p}}{2i}z_2\\
        0.1\left(\frac{p-\Bar{p}}{4i} + \frac{p+\Bar{p}}{4i}\right)\frac{z_1 + \Bar{z}_1}{2} + \Bar{p}\frac{z_1 - \Bar{z}_1}{2i} - \frac{3}{4}\left(\frac{1}{2} + \frac{1}{2i}\right)z_2
    \end{bmatrix}.
\end{equation}

\begin{equation}
    \partial_2^c g = \begin{bmatrix}
        \frac{2}{2i}\frac{p - \Bar{p}}{2i}\frac{z_1 + \Bar{z}_1}{2} + \left(\frac{10}{2}\frac{p+\Bar{p}}{2} + \frac{4}{2i}\frac{p-\Bar{p}}{2i}\right)\frac{z_1 - \Bar{z}_1}{2i} + 4\left(-\frac{1}{2i} - \frac{1}{2}\right)z_2\\
        -\frac{0.2}{2}\frac{p+\Bar{p}}{2}\frac{z_1 - \Bar{z}_1}{2i} + \frac{100}{2i}\frac{p-\Bar{p}}{2i}z_2\\
        0.1\left(\frac{p-\Bar{p}}{4i} - \frac{p+\Bar{p}}{4i}\right)\frac{z_1 + \Bar{z}_1}{2} + p\frac{z_1 - \Bar{z}_1}{2i} - \frac{3}{4}\left(\frac{1}{2} - \frac{1}{2i}\right)z_2
    \end{bmatrix}.
\end{equation}
The complex conjugates are applied element wise in these matrices. 
Using \autoref{eq: generalized adjoint equation}, we solve for $\frac{\partial z}{\partial p}$ and $\frac{\partial \Bar{z}}{\partial p}$. 
The partial derivatives of the cost function are the same as before, so $\nabla_z f = z$. 

\begin{sidewaysfigure}
    \centering
    \includegraphics[width=0.8\textheight]{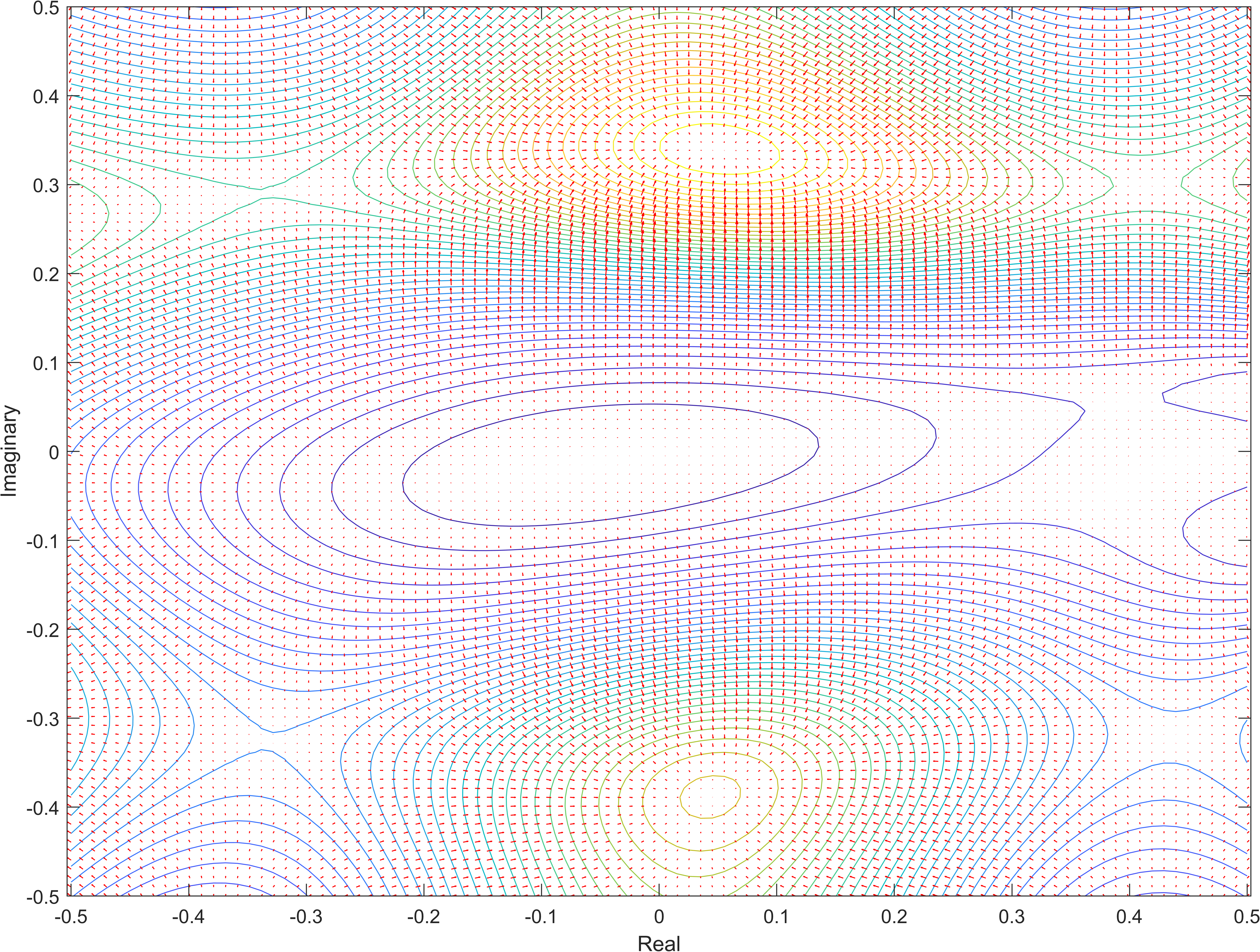}
    \caption{The gradient $\nabla_p f$ plotted as a vector field over the contour plot of $f$. The constraint written in \autoref{eq: non holomorphic constraint} is not holomorphic in both $p$ and $z$. }
    \label{fig: non-holomorphic contour plot}
\end{sidewaysfigure}

\begin{figure}
    \centering
    \includegraphics[width=0.75\linewidth]{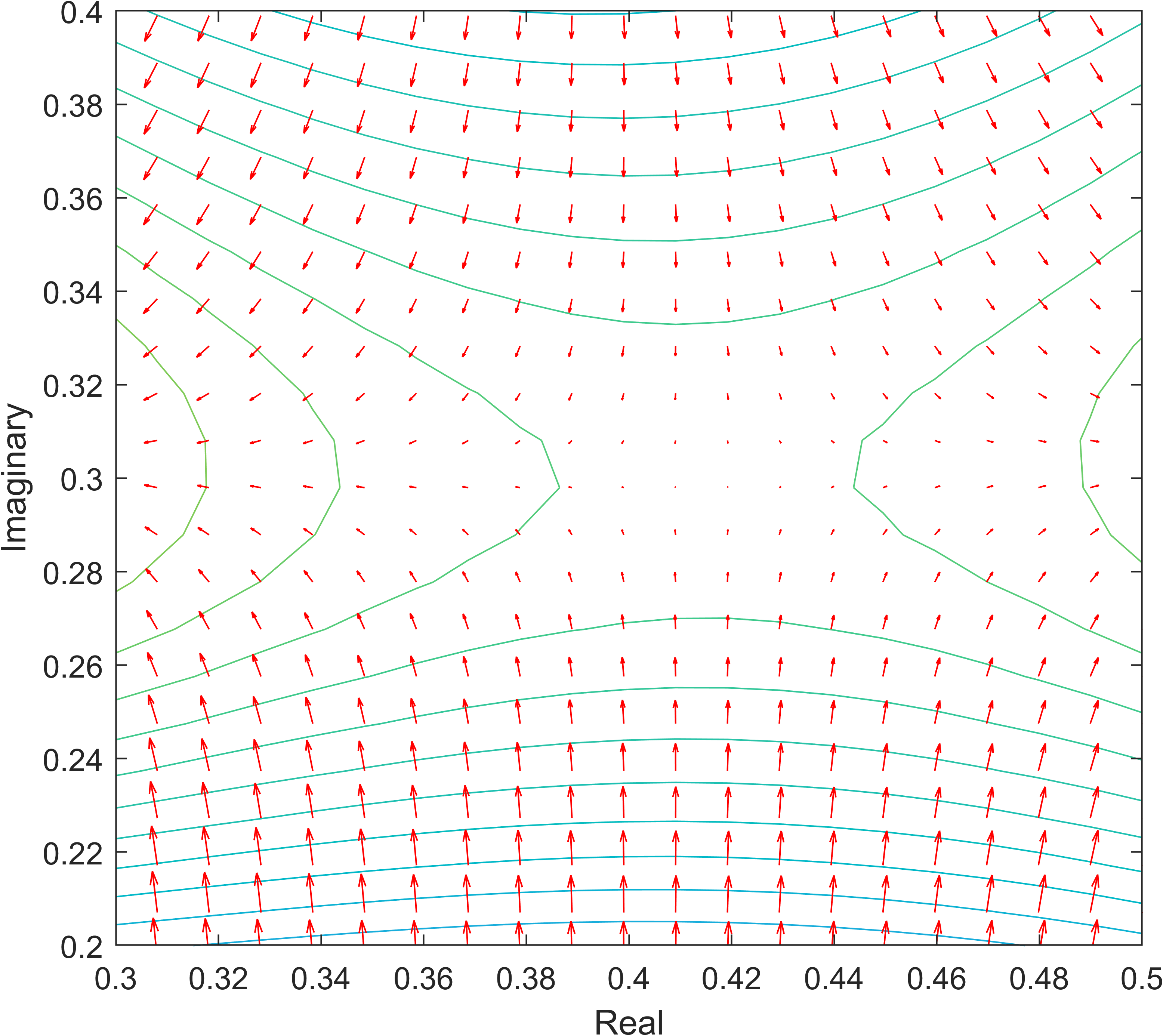}
    \caption{A zoomed in picture of \autoref{fig: non-holomorphic contour plot} showing the saddle point at $p \approx 0.4 + 0.3i$.}
    \label{fig: non-holomorphic contour plot saddle point}
\end{figure}

The gradients are computed using \autoref{eq: gradient for most general case} and plotted on a contour plot in \autoref{fig: non-holomorphic contour plot}. A zoomed in contour plot of the saddle point at $p \approx 0.4 + 0.3i$ is shown in \autoref{fig: non-holomorphic contour plot saddle point}.

\subsection{Example for 1-dimensional Helmholtz}

Consider the Helmholtz equation on the unit interval $[0, 1]$ with $f(x) = \sin(2\pi x)$
\begin{equation}\label{eq: Helmoltz}
    -\Delta u(x) - k^2 u = \sin(2\pi x)
\end{equation}
with Neumann boundary conditions $u'(0) = p_1, u'(1) = p_2$.
This can be solved using finite differences, which results in the PDE constraint being represented as a linear matrix 
\begin{equation}
    g(u, k, p_1, p_2) = A(k)u - b(p_1, p_2).
\end{equation}
Consider the case where 
\begin{equation}
    p_1 = p i, \quad p_2 = \bar{p}^3
\end{equation}
for some unknown $p \in \mathbb{C}$. This is not a realistic boundary condition, but is useful for showcasing the generalized adjoint method as $g$ is holomorphic in $u$ but not $p$. A more realistic boundary condition is $p_1, p_2\in \mathbb{C}$ but then the constraint would be holomorphic in $p = (p_1, p_2)\in \mathbb{C}^2$. 
The linear system is invertible, so we can define $u(k, p) = A(k)^{-1}b(p)$. 

Consider the inverse problem where $k^2 =4$ is fixed and we are trying to find an unknown $\Tilde{p}$. Let $\Tilde{u} = u(k, \Tilde{p})$. We choose the cost function to be 
\begin{equation}
    f(u, k, p) = |u(0) - \Tilde{u}(0)|^2 + |u(1) - \Tilde{u}(1)|^2.
\end{equation}

The partial derivatives required for the generalized adjoint method in this case are 
\begin{equation}
    \partial_1 g = A(k), \quad \partial_2 g = -\frac{\partial b(p)}{\partial p}, \quad \partial_2^c g = -\frac{\partial b(p)}{\partial \bar{p}}.
\end{equation}

Using finite differences with $N$ equidistance quadrature points and approximating the Neumann boundary condition with a second order forward/backwards difference
\begin{equation}
    u'(0) \approx -\frac{3}{2h}u(0) + \frac{2}{h}u(h) - \frac{1}{2h}u(2h), \quad u'(1) \approx \frac{3}{2h} u(1) - \frac{2}{h}u(1-h) + \frac{1}{2h}u(1 - 2h).
\end{equation}
Therefore, the partial derivatives of $b(p)$ in $p$ only affect the first and last element and we have 
\begin{equation*}
\begin{split}
    \frac{\partial b(p)}{\partial p} &= ie_1\\
    \frac{\partial b(p)}{\partial \bar{p}} &= 3\bar{p}^2 e_{N+2}
\end{split}
\end{equation*}
in which $e_k$ is the $k^{\text{th}}$ standard basis element for a $N+2$ dimensional vector space. 
The gradient $\nabla_p f$ can then be computed using \autoref{eq: adjoint for holomorphic in x but not p} and \autoref{eq: gradient for most general case}.

\begin{figure}
\begin{subfigure}[t]{0.45\textwidth}
    \centering
    \includegraphics[width=\textwidth]{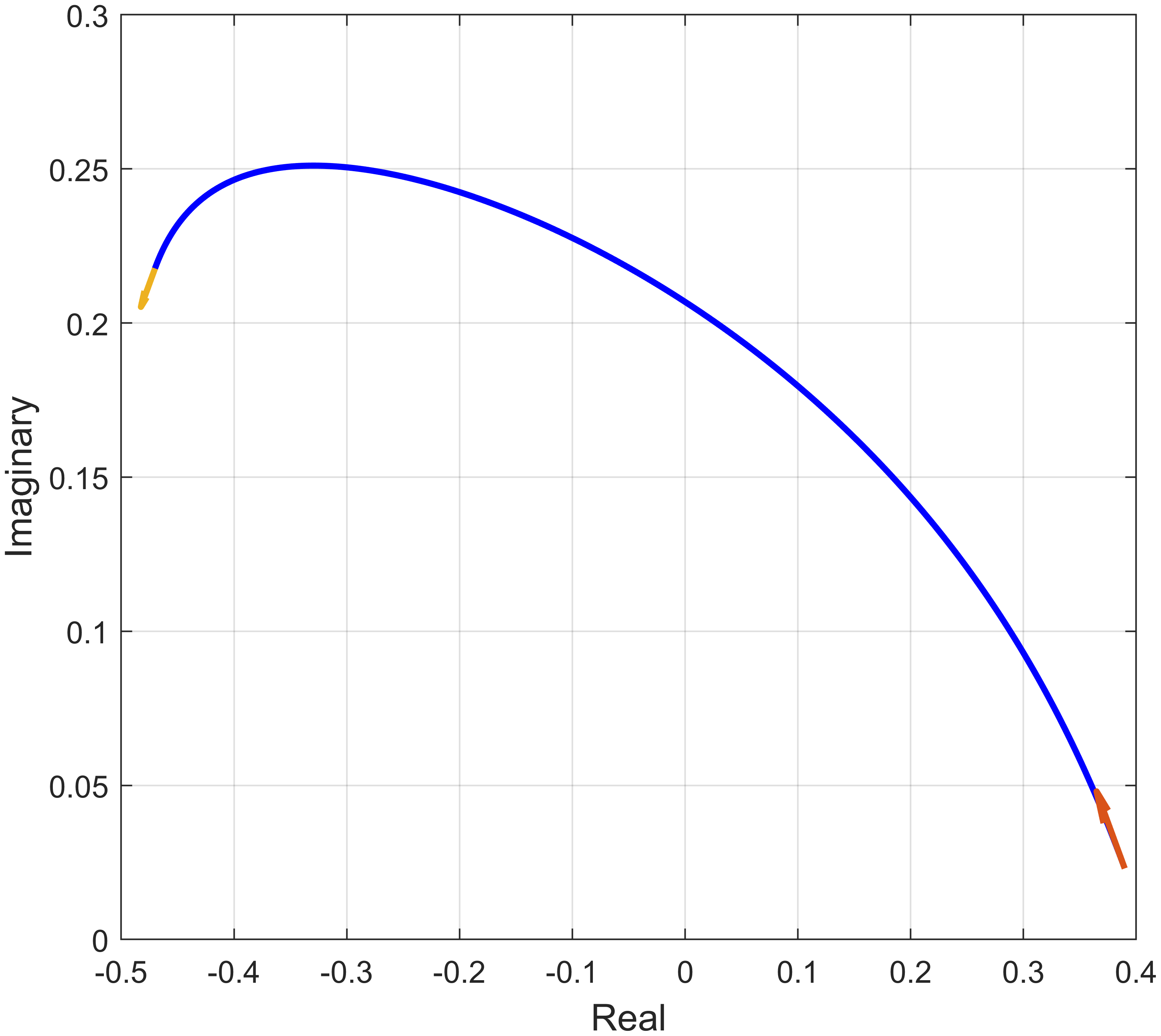}
    \caption{True solution of the 1D Helmholtz equation}
\end{subfigure}
\begin{subfigure}[t]{0.45\textwidth}
    \centering
    \includegraphics[width=\textwidth]{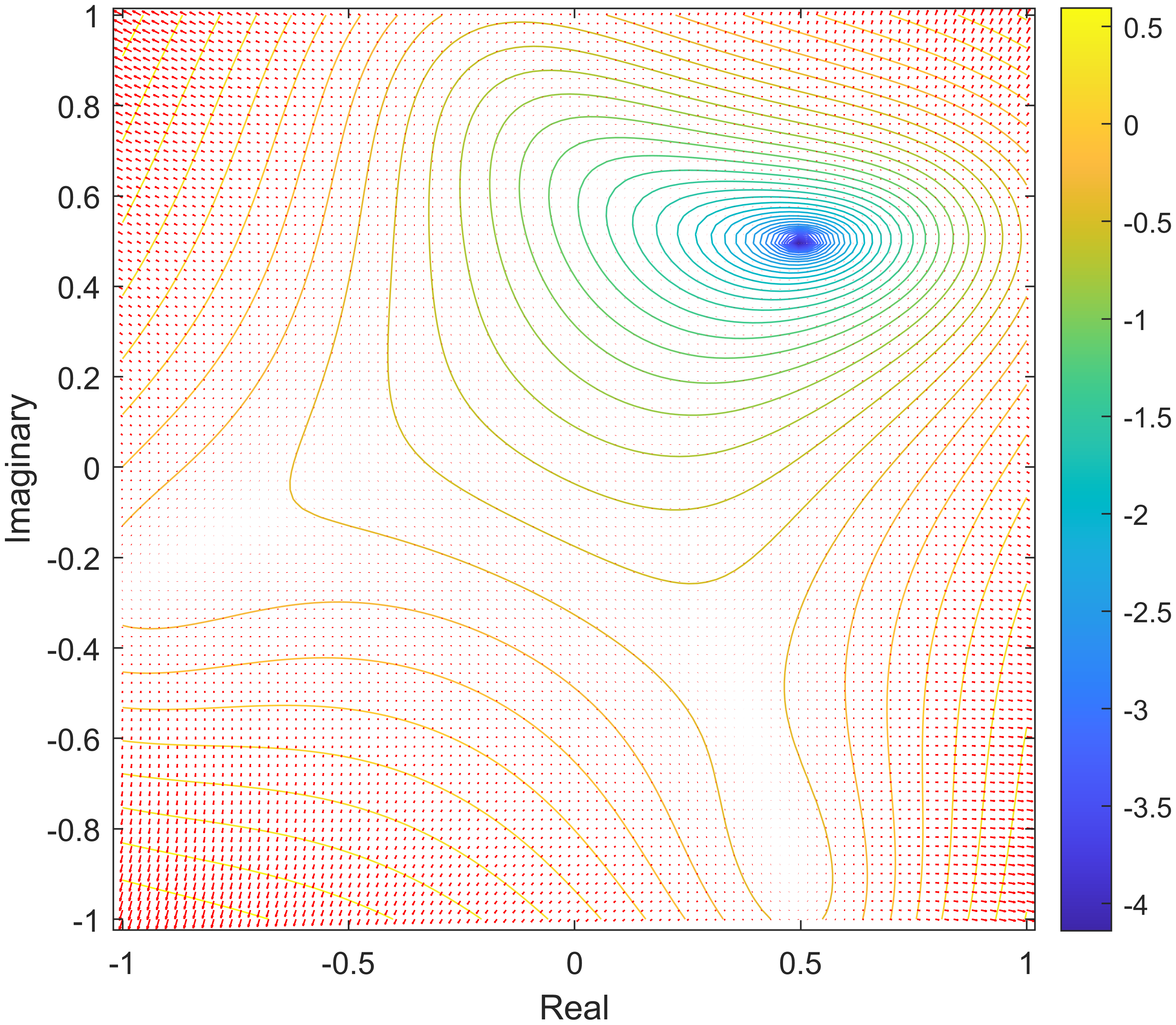}
    \caption{Contour plot of $\log_{10}(f(p))$.}
\end{subfigure}
    \caption{The true solution to the 1D Helmholtz equation with boundary condition with source term $\sin(2\pi x)$ and Neumann boundary condition $u'(0) = 0.5 - 0.5i$ and $u'(1) = -0.25 - 0.25i$ is shown in a). The Neumann boundary conditions are shown with vectors. The solution was obtained using finite difference with 1000 equi-spaced quadrature points. The contour plot of $\log_{10}(f(p))$ is shown in b) along with a vector field of the gradient. All gradient vectors are scaled equally so that the vectors do not overlap each other. The optimization plane is relatively flat close to the global minimum, which results in the vectors being very small in magnitude. Even when the optimizer gets very close to the global minimum, the slope is not that steep.  }
    \label{fig: Helmoltz contour and solution}
\end{figure}

With a true value of $p_{\text{true}} = 0.5 + 0.5i$ and an initial value of $p_0 = 0$, the inverse solver was solved using three standard optimization algorithms in MATLAB. The first test uses \texttt{fminunc} with a quasi-newton algorithm and we supply it with the gradient computed by the adjoint method. In the second test, we still use \texttt{fminunc}, but we do not give it a gradient which makes it default to finite differences. In the third test, we use the non-gradient based genetic algorithm. To make everything work in MATLAB, the domain of $p$ was written as $\mathbb{R}^2$, but the first step of the cost function is to rewrite it as a complex variable. This way, the cost function takes in real values and outputs a real value, but internally it uses complex variables. 
To not be convicted of the inverse crime \cite{wirgin2004inverse}, the true values of $\Tilde{u}(0)$ and $\Tilde{u}(1)$ were solved for using a quadrature of 1000 equi-spaced nodes, while the inverse problem was solved using 120 equi-spaced nodes. 
\autoref{fig: Helmoltz contour and solution} shows how the optimization landscape and the true solution to the Helmholtz equation in \autoref{eq: Helmoltz} behave like. 
\autoref{tab: Helmholtz inverse results} shows that with the generalized adjoint method, the number of function calls is greatly reduced. For larger problems, or more quadrature points, a reduction in the number of function calls greatly speeds up the optimization algorithm. 
\begin{table}[]
    \centering
    \begin{tabular}{|c|c|c|c|}
    \hline
        Algorithm & \makecell{Number of \\ function calls} & Final $p$ &Final cost \\
        \hline
        \makecell{\texttt{fminunc} with\\generalized adjoint method} & 13 & $0.50002 + 0.49987i$ & $4.11\cdot 10^{-10}$\\
        \hline
        \makecell{\texttt{fminunc} with \\ finite differences in $p$} & 70 & $0.50000 + 0.49987i$ & $4.09\cdot 10^{-10}$\\
        \hline
        genetic algorithm & 5615 & $0.50004 + 0.49989i$ & $8.93\cdot 10^{-10}$\\
        \hline
    \end{tabular}
    \caption{Results of solving the Helmholtz 1D inverse problem using the generalizd adjoint method, finite differences, and the non-gradient based genetic algorithm. The global minimum has a cost of $0$.}
    \label{tab: Helmholtz inverse results}
\end{table}

Tests with small $k^2 \ll 1$ (high frequency) result in an ill-posed problem \cite{babuska1997pollution}. In this regime, \texttt{fminunc} fails with both the adjoint method and finite differences. Only the genetic algorithm is able to find the true global minimum. 


\section{Conclusion}\label{sec: conclusion}

Many optimization problems represent certain parameters and signals using complex variables. This results in a better model of the real world problem, but comes at the price of introducing non-holomorphic functions into the optimization problem. Standard gradient based optimization methods no longer apply, and representing $\mathbb{C}^n$ as $\mathbb{R}^{2n}$ is not recommended. In the specific case of equality constraint optimization, we introduce the generalized adjoint method to deal with non-holomorphic cost functions and constraint functions. 
Using CR-calculus, the generalized adjoint method gives equations to calculate the gradient of non-holomorphic functions. Combining this with differential geometry or an equivalent Lagrangian formulation, the optimizer is able to accurately determine the direction of steepest ascent/descent while enforcing the non-holomorphic constraint. The equations also accurately capture the different cases of whether the domain over which we optimize is real or complex. 

Further extensions of this generalized adjoint method to compute Hessians and higher order derivatives can be done with the Lagrangian formulation in a method similar to \cite{antil2018brief}, but are much more complicated and hence not included in this paper. 

Another potential concern is the difficulty of evaluating certain partial derivatives numerically. In practice, finite difference methods can be used at the cost of high error. Automatic differentiation is another potential technique that can be used, but may not work for certain types of code \cite{rall1981automatic}. In most cases, computing the partials of the constraint function is easier than computing $\frac{\partial x}{\partial p}$ and $\frac{\partial x}{\partial \bar{p}}$. In the case of PDE constraints, the difficulty of partial derivatives will be different using a continuous or discrete formulation, and it is up to the user to decide which formulation is used. 





\section*{Declaration of competing interest}
The authors have no competing interests to declare that are relevant to the content of this article.




\bibliographystyle{elsarticle-num} 
\bibliography{bibliography}



\end{document}